\documentclass[12pt]{article}
\usepackage{amssymb}
\usepackage{amsmath}
\usepackage[active]{srcltx}
\newtheorem{theorem}{Theorem}[section]
\newtheorem{lemma}[theorem]{Lemma}

\newtheorem{conjecture}[theorem]{Conjecture}
\newtheorem{proposition}[theorem]{Proposition}
\newcommand{\qed}{\hfill $\square$}
\newcommand{\comments}[1]{}

\DeclareMathOperator* \ave {ave}
\newcommand{\hb}{\bar{h}}
\newcommand{\A}{\mathcal{A}}
\newcommand{\B}{\mathcal{B}}
\newcommand{\C}{\mathcal{C}}
\newcommand{\D}{\mathcal{D}}
\newcommand{\E}{{\rm E}}
\newcommand{\F}{\mathcal{F}}

\newcommand{\Ha}{\mathcal{H}}
\newcommand{\J}{\mathcal{J}}
\newcommand{\K}{\mathcal{K}}
\newcommand{\La}{{\rm La}}
\newcommand{\M}{\mathcal{M}}
\newcommand{\N}{\mathcal{N}}
\newcommand{\Oh}{\mathcal{O}}
\newcommand{\Pa}{\mathcal{P}}
\newcommand{\V}{\mathcal{V}}

\newcommand{\nchn}{\binom{n}{\lfloor \frac{n}{2}\rfloor}}
\newcommand{\mchm}{\binom{m}{\lfloor \frac{m}{2}\rfloor}}

\newcommand{\lanp}{\La(n,P)}

\begin{document}
\title{Diamond-free Families}
\author{Jerrold R. Griggs\thanks{
Department of Mathematics, University of South Carolina, Columbia,
SC 29208  USA ({\tt griggs@math.sc.edu}).} \and Wei-Tian
Li\thanks{Department of Mathematics, University of South Carolina,
Columbia SC 29208 USA ({\tt li37@mailbox.sc.edu}).} \and Linyuan
Lu
\thanks{Department of Mathematics,
University of South Carolina, Columbia, SC 29208  USA ({\tt
lu@math.sc.edu}). This author was supported in part by NSF grant
DMS 1000475. } }

\date{August 28, 2011\\}

\maketitle

\begin{abstract}
Given a finite poset $P$, we consider the largest size $\lanp$ of
a family of subsets of $[n]:=\{1,\ldots,n\}$ that contains no
(weak) subposet $P$.
This problem has been studied intensively in recent years, and it
is conjectured that $\pi(P):= \lim_{n\rightarrow\infty}
\lanp/\nchn$ exists for general posets $P$, and, moreover, it is
an integer. For $k\ge2$ let $\D_k$ denote the $k$-diamond poset
$\{A< B_1,\ldots,B_k < C\}$. We study the average number of times
a random full chain meets a $P$-free family, called the Lubell
function, and use it for $P=\D_k$ to determine $\pi(\D_k)$ for
infinitely many values $k$. A stubborn open problem is to show
that $\pi(\D_2)=2$; here we make progress by proving
$\pi(\D_2)\le 2\frac{3}{11}$ (if it exists).
\end{abstract}

\section{Introduction}
We are interested in how large a family of subsets of the $n$-set
$[n]:=\{1,\ldots,n\}$ there is that avoids a given (weak) subposet
$P$. The foundational result of this sort, Sperner's Theorem from
1928 ~\cite{Spe}, solves this problem for families that contain no
two-element chain (that is, for antichains), determining that the
maximum size is $\nchn$.  For other excluded subposets, it is
interesting to compare the maximum size of a $P$-free family to
$\nchn$.

We give background and our new results for this study in the next
section. One small forbidden poset that continues to stymie all
interested researchers is the diamond poset on four elements.  We
present a better new upper bound on the size of diamond-free
families.  For $k$-diamond-free families for general $k$, we
provide bounds that, surprisingly, turn out to be best-possible
for infinitely many values of $k$.

In Section 3 we introduce our method for this subject, the Lubell
function of a family, which gives the average number of times a
random full chain meets the family.  The Lubell function yields an
upper bound on the size of a family.  For diamond-free families,
we observe that the maximum possible Lubell function decreases
with $n$, and a calculation gives our bound.  For excluding
$k$-diamonds, our new idea is to partition the set of full chains,
obtaining bounds on each block of the partition.  The Lubell and
full chain partition methods hold promise for other families of
forbidden subposets.

Section 4 contains our detailed proofs, except for the long proof
of Theorem~\ref{d2better}, which is given in Section~5.  The paper
concludes with our ideas for advancing the project.

\section{Background and Main Results}

For posets $P=(P,\le)$ and $P'=(P',\le')$, we say $P'$ is a {\em
weak subposet of $P$} if there exists an injection $f\colon P'\to P$
that preserves the partial ordering, meaning that whenever $u\le'
v$ in $P'$, we have $f(u)\le f(v)$ in $P$~\cite{Sta}. Throughout the
paper, when we say subposet, we mean weak subposet. We say the {\em
height\/} $h(P)$ of poset $P$ is the maximum size of any chain in
$P$.

Let the Boolean lattice $\B_n$ denote the poset
$(2^{[n]},\subseteq)$.  We consider collections $\F\subseteq
2^{[n]}$.  Then $\F$ can be viewed as a subposet of $\B_n$. If
$\F$ contains no subposet $P$, we say $\F$ is {\em $P$-free.}  We
are interested in determining the largest size of a $P$-free
family of subsets of $[n]$, denoted $\La(n,P)$.

In this notation, Sperner's Theorem~\cite{Spe} gives that
$\La(n,\Pa_2)=\nchn$, where $\Pa_k$ denotes the path poset on $k$
points, usually called a chain of size $k$. Moreover, Sperner
determined that the largest antichains in $\B_n$ are the middle
level, $\binom{[n]}{{n/2}}$ (for even $n$) and either of the two
middle levels, $\binom{[n]}{(n-1)/2}$ or $\binom{[n]}{(n+1)/2}$
(for odd $n$), where for a set $S$, $\binom{S}{i}$ denotes the
collection of $i$-subsets of $S$. More generally, Erd\H os solved
the case of $\Pa_k$-free families. Let us denote by $\Sigma(n,k)$
the sum of the $k$ middle binomial coefficients in $n$, and let
$\B(n,k)$ denote the collection of subsets of $[n]$ of the $k$
middle sizes, that is, the  sizes $\lfloor (n-k+1)/2 \rfloor,
\ldots, \lfloor (n+k-1)/2 \rfloor$ or else the sizes $\lceil
(n-k+1)/2 \rceil, \ldots, \lceil (n+k-1)/2 \rceil$. So there are
either one or two possible families $\B(n,k)$, depending on the
parities of $n$ and $k$, and regardless, $|\B(n,k)|=\Sigma(n,k)$.
Then we have

\begin{theorem}~\cite{Erd}\label{Erd}
For $n\ge k-1\ge 1$, $\La(n,\Pa_k)=\Sigma(n,k-1)$.  Moreover, the
$\Pa_k$-free families of maximum size in $\B_n$ are given by
$\B(n,k)$.
\end{theorem}

It follows that for fixed $k$, $\La(n,\Pa_k)\sim (k-1)\nchn$ as
$n\rightarrow\infty$.  Katona and his collaborators promoted the
problem of investigating $\La(n,P)$ for other posets $P$,
especially its asymptotic behavior for large $n$.  Consider the
$r$-fork poset $\V_r$, which has elements $A<B_1,\ldots,B_r$,
$r\ge2$.  In 1981 he and Tarj\'an~\cite{KatTar} obtained bounds
on $\La(n,\V_2)$ that he and DeBonis~\cite{DebKat}  extended in
2007 to general $\V_r$, $r\ge2$, proving that
\[
\left(1+\frac{r-1}{n} + \Omega\left(\frac{1}{n^2}\right)\right)\nchn
\le\La(n,\V_r)\le \left(1+ 2\frac{r-1}{n} +
O\left(\frac{1}{n^2}\right)\right)\nchn.
\]
While the lower bound is
strictly greater than $\nchn$, we see that $\La(n,\V_r)\sim \nchn$.
Earlier, Thanh~\cite{Tha} had investigated the more general class
of broom-like posets. Griggs and Lu~\cite{GriLu} studied the even
more general class of baton posets.  These posets mentioned so far
have Hasse diagrams that are trees.

In~\cite{DebKatSwa} it is shown that for the butterfly poset $\B$,
with elements $A,B$ both less than $C,D$, one can give an exact
answer, $\La(n,\B)=\Sigma(n,2)$, for $n\ge3$, which is asymptotic
to $2\nchn$.  More generally, for any $s,t\ge 2$, the complete
bipartite poset $\K_{s,t}$ with elements $A_1,\ldots,A_s$ all less
than $B_1,\ldots,B_t$, satisfies $\La(n,\K_{s,t})\sim 2\nchn$
\cite{DebKat}. The $\N$-poset, with elements $A, B, C, D$ such
that $A<B, C<B, C<D$, is intermediate between $\V_2$ and $\B$.  It
is shown in~\cite{GriKat} that $\La(n,\N) \sim \nchn$.

Based on the examples for which $\lanp$ was known, Griggs and
Lu~\cite{GriLu} proposed the conjecture that was certainly
apparent to Katona {\em et al.}:

\begin{conjecture}\label{conj}
For every finite poset $P$, the limit
$\pi(P):=\lim_{n\rightarrow\infty} \frac{\lanp}{\nchn}$ exists and
is integer.
\end{conjecture}

All of the examples above agree with the conjecture, and Griggs
and Lu verified it for additional examples, including tree posets
of height 2.  For the crown $\Oh_{2k}$, $k\ge2$, which is the
poset of height $2$ that is a cycle of length $2k$ as an
undirected graph, they extended the butterfly result above and
proved that $\La(n,\Oh_{2k})\sim \nchn$ for all even $k\ge4$. For
odd $k\ge3$, it remains a daunting problem to determine the
asymptotic behavior of $\La(n,\Oh_{2k})$. At least, Griggs and Lu
can show $\La(n,\Oh_{2k})/\nchn$ is asymptotically at most
$1+\frac{1}{\sqrt2}$, which is less than 2.

When Griggs lectured on this work on forbidden subposets in 2008,
Mike Saks and Peter Winkler observed a pattern in all of the
examples where $\pi(P)$ was determined, which we describe as
follows.  For poset $P$ define $e(P)$ to be the maximum $m$ such
that for all $n$, the union of the $m$ middle levels $\B(n,m)$
does not contain $P$ as a subposet. Their observation was
$\pi(P)=e(P)$.  For instance, the middle two levels $\B(n,2)$
contain no butterfly $\B$, since no two sets of the same size $k$
contain the same two subsets of size $k-1$.  One gets that
$e(\B)=2$, which is $\pi(\B)$. In general, it is clear that when
it exists, $\pi(P)$ must be at least $e(P)$.

Impressive progress in the development of the theory is the result
of Bukh~\cite{Buk} that for any tree poset $T$ (meaning that the
Hasse diagram is a tree), $\pi(T)=e(T)$, so that the conjecture
(and observation) are satisfied.  It is easily verified that
$e(T)=h(T)-1$.

Is there a connection for general $P$ between $\pi(P)$ and the
height $h(P)$? The result of DeBonis and Katona for complete
bipartite posets $\K_{s,t}$ implies that for any poset $P$ of
height 2, $\pi(P)\le2$, when it exists. However, there is no such
bound for taller posets, as observed by Jiang and Lu (see
~\cite{GriLu}). Let the $k$-diamond poset $\D_k$, $k\ge2$, consist
of $k+2$ elements $A<B_1,\ldots,B_k<C$. Then $h(\D_k)=3$ for
general $k$, while for $k=2^r-1$, the middle $r+1$ levels
$\B(n,r+1)$ cannot contain $\D_k$, since an interval in
$\B(n,r+1)$ with an element in the lowest level and an element in
the highest level has at most $2^r$ elements (a subposet $\B_r$),
and so at most $2^r-2$ elements in the middle. Hence it is
$\D_k$-free.

The diamond $\D_2$ is the most challenging poset on at most four
elements in this theory.  (It is also the Boolean lattice $\B_2$.)
It is easily seen that $e(\D_2)=2$.  On the other hand, it is a
subposet of the path $\Pa_4$.  So, if $\pi(\D_2)$ exists (which
has still not been shown), it would have to be in $[2,3]$; Its
conjectured value is 2.

As an illustration of the Lubell function method introduced in the
next section, a short application is given that reduces the upper
bound on $\pi(\D_2)$ from 3 to $2.5$.  A refinement of the Lubell
function method, which involves partitioning the set of full
chains in an appropriate way, gives our first improvement on the
$2.5$ bound:

\begin{proposition}\label{d2}
For all sufficiently large $n$, $\La(n,\D_2)/\nchn < 2.296$.
\end{proposition}

We display this bound, not our best one, since its proof is
simpler than our best bound, and since its proof gives us further
insight into the Lubell function for $\D_2$.  Some time after we
had announced our bound above, Axenovich, Manske, and
Martin~\cite{AxeManMar} came up with a new approach which improves
the upper bound to $2.283$.  Now using our methods with a much
more careful analysis of diamond-free families for $n\le12$, we
can provide a further slight improvement, which is the best-known
upper bound:

\begin{theorem}\label{d2better}
For all sufficiently large $n$, $\La(n,\D_2)/\nchn < 2\frac{3}{11}+o_n(1)<2.273$.
Consequently, if it exists, $\pi(\D_2)\in [2,2\frac{3}{11}]$.
\end{theorem}

Because this new bound requires considerably more care, its proof
is given in its own section following the proofs of our other
results.  We shall see diamond-free families in the proof for
which the Lubell function method cannot improve the upper bound on
$\pi(\D_2)$ below $2.25$.  Therefore, new ideas are required to
bring the upper bound down to the conjectured value of 2.
Likewise, it appears that the methods of~\cite{AxeManMar} cannot
move below $2.25$.  See the final section of the paper for more
discussion of how we can do better.

Given the great effort that has gone into improving the upper
bound on $\pi(\D_2)$, it is then quite surprising that we can
solve the $\pi$ problem for many of the general diamonds $\D_k$
with $k>2$. This can be regarded as our main result.

\begin{theorem}\label{dk}
Let $k\ge2$, and define $m:=\lceil \log_2(k+2)\rceil$.

(1)  If $k\in [2^{m-1}-1, 2^m-\mchm-1]$, then
\[
\La(n,\D_k)=\Sigma(n,m).
\]
Hence, $\pi(\D_k)=e(\D_k)=m$.
Moreover, if $\F$ attains the bound $\La(n,\D_k)$, then
$\F=\B(n,m)$.

(2)  If $k\in [2^m-\mchm, 2^m-2]$, then,
\[
\Sigma(n,m)\le
\La(n,\D_k)\le \left(m+1-{\frac{2^m-k-1}{\mchm}}\right) \nchn.
\]
Hence, if $\pi(\D_k)$ exists, then
\[
m=e(\D_k)\le \pi(\D_k)\le
\left(m+1-{\frac{2^m-k-1}{\mchm}}\right)<m+1.
\]
\end{theorem}

For  $\pi(\D_2)$ this new theorem gives an upper bound of $2.5$,
not as good as the theorem before.  However, this new result
determines $\pi(\D_k)$ for ``most" values of $k$, in that for the
$2^{m-1}$ values of $k$ in the range $[2^{m-1}-1,2^m-2]$, case (1)
applies to all but $\mchm\sim C2^m/m^{1/2}$ of them.  Moreover, we
are able to give $\La(n,\D_k)$ exactly, not just asymptotically
for large $n$, for such values of $k$.

The poset $\D_k$ can be viewed as the ``suspension" of an
independent set of size $k$, where we mean that a maximum and a
minimum element are added to it.  We can consider a more general
suspension of disjoint paths (chains).  For $k\ge1$ let $l_1\ge
\cdots\ge l_k\ge3$, and define the {\em harp poset\/}
$\Ha(l_1,\ldots,l_k)$ to consist of paths
$\Pa_{l_1},\ldots,\Pa_{l_k}$ with their top elements identified
and their bottom elements identified.  For instance, in this
notation we have $\D_k$ is the harp $\Ha(3,\ldots,3)$ where there
are $k$ 3's.

\begin{theorem}\label{harp}
If $l_1>\cdots>l_k\ge3$, then
\[
\La(n,\Ha(l_1,\ldots,l_k))=\Sigma(n,l_1-1).
\]
Hence, for such harps, $\pi=e=l_1-1$. Moreover, for such harps, if
$\F$ is a harp-free family of subsets of $[n]$ of maximum size,
then $\F$ is $\B(n,l_1-1)$.
\end{theorem}

The theorem above only determines $\pi(\Ha)$ for harps $\Ha$ that
have strictly decreasing path lengths.  However, for the general
case in which path lengths can be equal there is no bound
independent of $k$, since we have seen that for $\D_k$, which is a
harp, $\pi(\D_k)$ is arbitrarily large as $k$ grows.  It is then
remarkable that we can completely solve the problem of maximizing
$\La(n,\Ha)$ for harps with distinct path lengths. Another novel
aspect of this result is that for $k\ge2$ the harps it concerns
are not ranked posets.

\section{The Lubell Function}

For now let us fix some family $\F\subseteq 2^{[n]}$.  Let
$\C:=\C_n$ denote the collection of all $n!$ full (maximal) chains
$\emptyset \subset \{i_1\} \subset \{i_1,i_2\} \subset \cdots
\subset [n]$ in the Boolean lattice $\B_n$. A method used by
Katona {\em et al.\/} involves counting the number of full chains
that meet $\F$.  Here we collect information about the average
number of times chains $C\in \C$ meet $\F$, which can be used to
give an upper bound on $|\F|$. Recall that the {\em height\/} of
$\F$, viewed as a poset, is
\[
h(\F):= \max_{C\in\C} |\F\cap C|.
\]
We consider what we call the {\em Lubell function\/} of $\F$,
which is
\[
\hb(\F)=\hb_n(\F):=\ave_{C\in\C} |\F\cap C|.
\]
This is the expected value $\E(|\F\cap C|)$ over a random full
chain $C$ in $\B_n$.  Then $\hb(\F)$ is essentially the function
of $\F$ at the heart of Lubell's elegant proof of Sperner's Theorem
(\cite{Lub}, cf. \cite{GreKle}) with the observation.

\begin{lemma}\label{lub}
Let $\F$ be a collection of subsets of $[n]$.  Then $\hb(\F) =
\sum_{F\in \F} 1/\binom{n}{|F|}$.
\end{lemma}
{\bf Proof}: We have that $\hb(\F)= \E(|\F\cap C|)$, where $C$ is
picked at random from $\C$.  This expected value is, in turn, the
sum over $F\in \F$ of the probability that a random $C$ contains
$F$.  Since $C$ meets the $\binom{n}{k}$ subsets of cardinality
$k$ with equal probability, it means that each set $F$ contributes
$1/\binom{n}{|F|}$ to the sum. \qed

\medskip
Lubell's proof uses the simple facts that $|\A\cap C|\le 1$ for
any antichain $\A$ and that $\binom{n}{k}$ is maximized by taking
$k=\lfloor \frac{n}{2}\rfloor$, to derive Sperner's Theorem that
$|\A|\le \nchn$. By similar reasoning for general families $\F$ we
obtain a general upper bound.

\begin{lemma}\label{lubbound}
Let $\F$ be a collection of subsets of $[n]$.  If $\hb(\F)\le m$,
for real number $m>0$, then $|\F|\le m\nchn$.  Moreover, if $m$ is
an integer, then $|\F|\le\Sigma(n,m)$, and equality holds if and
only if $\F=\B(n,m)$ (when $n+m$ is odd), or if $\F=\B(n,m-1)$
together with any $\binom{n}{(n-m)/2}$ subsets of sizes $(n-m)/2$
or $(n+m)/2$ (when $n+m$ is even).
\end{lemma}
{\bf Proof}: We use the symmetry and strict unimodality of the
sequence of binomial coefficients $\binom{n}{k}$, $0\le k\le n$.
If $\hb(\F)\le m$, then $|\F|=\sum_{A\in\F}1\le
\sum_{A\in\F}\nchn/\binom{n}{|A|}\le m\nchn$. Now assume $m>0$ is
an integer.  We construct a family $\F$ of maximum size, subject
to $\hb(\F)\le m$, by selecting subsets $A$ that contribute the
least to $\hb(\F)$, which means that we minimize
$1/\binom{n}{|A|}$. Essentially, we are solving the linear
program of maximizing $\sum_i x_i \binom{n}{i}$ subject to
$\sum_i x_i\le m$, $0\le x_i\le 1$ for all $i$.  We maximize $|\F|$ by
selecting $\F$ to be the $m$ middle levels, $\B(n,m)$.

Further, if $|\F|=\Sigma(n,m)$, it must be that $\F$ is $\B(n,m)$
when $n+m$ is odd. If $n+m$ is even, the subsets of sizes
$(n-m)/2$ and $(n+m)/2$ will tie for the $m$-th largest size, and
we can freely choose any $\binom{n}{(n-m)/2}$ subsets of the two
sizes so that $|\F|=\Sigma(n,m)$. \qed

\medskip
We see that upper bounds on the average intersection size $|\F\cap
C|$ lead to upper bounds on the ratio of particular interest in
this paper, $|\F|/\nchn$.   Hence, we get upper bounds on
$\pi(P)$, when it exists, from upper bounds on $\hb(\F)$ for
$P$-free families $\F$.

To illustrate how this can be useful, we now give a short proof
that, if it exists, $\pi(\D_2)\le 2.5$.  Consider a diamond-free
family $\F\subseteq 2^{[n]}$.  No full chain $C\in\C_n$ meets $\F$
four times, or else $\F$ contains $\Pa_4$, which has $\D_2$ as a
subposet.  If no chain meets $\F$ three times, we immediately get
$\hb(\F)\le h(\F)\le2$.  Else, consider any three elements of $F$
$X\subset Y\subset Z$, and let $Y'$ be any set not equal to $Y$
such that $X\subset Y'\subset Z$.  Let $\sigma$ be a permutation
of $[n]$ that fixes $X$ and $Z$ and sends $Y$ to $Y'$.  Then
$\sigma$ sends full chains $C\in \C_n$ through $X,Y,Z$ to full
chains $C'\in\C$ through $X,Y',Z$.  These chains $C'$ meet $\F$
only twice, as $\F$ is diamond-free.  These chains $C'$ are
distinct.  We find then that $\hb(\F)=\E(|\F\cap C|)\le 2.5$.

Unfortunately, the behavior of $\hb(\F)$ does not match that of
$|\F|/\nchn$ asymptotically--there can be a gap. We shall see
examples of this for diamond-free families. Nonetheless, in many
cases we can obtain $\pi(P)$ from $\hb(\F)$. Besides that, it is
interesting in its own right to maximize $\hb(\F)$ for $P$-free
families $\F$, though obtaining a good bound on $\hb(\F)$ can be
difficult. We have discovered that a ``partition method'' can be
fruitful.

Specifically, we partition the set $\C$ of full chains into blocks
$\C(i)$ and then, for each $i$ separately, we bound the average
size $|\F\cap C|$ over full chains $C\in \C(i)$. The principle is
that the average size $|\F\cap C|$ over all full chains $C$ is at
most the maximum over $i$ of the average over block $\C(i)$. An
analogy to baseball is helpful for some readers: A hitter's
average over a whole season is never more than his maximum monthly
average over the months in the season.

We illustrate the partition method by sketching a derivation of
$\La(n,\B)$. Let $\F$ be a butterfly-free family of subsets of
$[n]$, $n\ge3$. One can check that if $\F$ contains $\emptyset$ or
$[n]$, then $|\F|<\Sigma(n,2)$ (although, one may have
$\hb(\F)>2$).  Else, suppose $\emptyset,[n]\notin \F$. We show
$\hb(\F)\le2$.  Define the collection $\M$ of subsets $M\in \F$ of
$[n]$ for which there exists a chain $C\in \C$ passing through
$A,M,B\in \F$ with $A\subset M\subset B$.  Notice that since $\F$
contains no butterfly $\B$, it contains no $\Pa_4$, and so the
collection $\M$ is an antichain.  Now partition the set of full
chains $\C$ as follows:  For $M\in \M$, $\C_M$ consists of all
full chains meeting $M$, while $\C_{\emptyset}$ contains all full
chains that do not meet $\M$.

By definition of $\M$, no chain in $\C_{\emptyset}$ meets $\F$
three times, and so $\ave_{C\in \C_{\emptyset}} |\F\cap C|\le 2$.
For $M\in \M$, similar to the argument above for $\D_2$-free
families, for any chain $C\in\C_M$ meeting $\F$ three times, it
must meet $\F$ in $A,M,B$, and there is a corresponding chain in
$\C_M$ meeting $\F$ only at $M$ and avoiding $A,B$, so that
$\ave_{C\in \C_M} |\F\cap C|\le 2$.  Hence, we have partitioned
$\C$ into blocks such that $\F$ meets chains in each block at most
twice, on average, and hence at most twice, on average, over all
of $\C$.  Thus, $\hb(\F)\le 2$, and it follows that
$\La(n,\B)=\Sigma(n,2)$, since $\B(n,2)$ is butterfly-free.

Regarding the extremal butterfly-free families as far as achieving
$\La(n,\B)$, Lemma~\ref{lubbound} above applies.  In fact, it is
known that $\F$ must be $\B(n,2)$ for $n\ge5$, though it is not
true for $n=4$:  Consider $\F$ consisting of
$\{\{1\},\{2\},\{1,3,4\},\{2,3,4\}\}$ and all six 2-subsets.

However, in some cases we can show that $\La(n,P)$ is attained
only by $\F=\B(n,k)$:

\begin{lemma}\label{exfam}
Suppose that for poset $P$, $e(P)=m$, an integer.  Suppose that
for all $n$, all $P$-free families $\F\subseteq 2^{[n]}$ satisfy
$\hb(\F)\le m$.  Then for all $n$, $\lanp=\Sigma(n,m)$, and if
$\F$ is an extremal family, then $\F=\B(n,m)$.
\end{lemma}
{\bf Proof}: Let $\F$ be a $P$-free family with size $\La(n,P)$.
According to Lemma \ref{lubbound}, $\La(n,P)=\Sigma(n,m)$, since
$\hb(\F)\le m$.  Further, $\F=\B(n,m)$ when $n+m$ is odd.  Hence,
suppose $n+m$ is even, so that $\F=\B(n,m-1)$ together with any
$\binom{n}{(n-m)/2}$ subsets of sizes $(n-m)/2$ or $(n+m)/2$.
Suppose for contradiction that $\F$ contains subsets of both sizes
$(n-m)/2$ and $(n+m)/2$.  By the natural generalization of
Sperner's proof of Sperner's Theorem (or by using the normalized
matching property on the ranks $\binom{[n]}{(n-m)/2}$ and
$\binom{[n]}{(n+m)/2}$)~\cite{GreKle}, we can find subsets $A,B\in
\F$ with $A\subset B$, $|A|=(n-m)/2$ and $|B|=(n+m)/2$.  Then
the interval $[A,B]\subset \F$.  This interval is a Boolean
lattice, $\B_m$.  However, $\hb_m(\B_m)=m+1>m$, so that by
hypothesis, $\B_m$ must contain subposet $P$, which contradicts
our assumption that $\F$ is $P$-free.  Hence, $\F$ only contains
one of the two sizes $(n-m)/2$ and $(n+m)/2$.\qed

\medskip
 We saw that the conclusion of the lemma above fails for
the butterfly $P=\B$ (but only for small $n$).  The reason we
could not apply this lemma to $P=\B$ is that the hypothesis fails
for $n=2$: The full Boolean lattice $\B_2$, which has one more
element than $\B(2,2)$, is butterfly-free.

\section{Proofs of Results \ref{d2}, \ref{dk}, \ref{harp}}

We now illustrate our partition method to bring the bound for
$\D_2$-free families below $2.5$.  Our best bound is derived in
the next section.
\medskip

\noindent{\bf Proof of Proposition~\ref{d2}:}

Let $\F$ be a $\D_2$-free family of subsets of $[n]$ with maximum
Lubell function value $\hb(\F)$, and let $d_n$ denote this value.
We claim that $d_n$ is nonincreasing for $n\ge2$.  By easy direct
case study we get that $d_2=2.5$ and $d_3=d_4=7/3\approx2.33$.

For $n\ge 5$, if both $\emptyset$ and $[n]$ are in $\F$, then we
have only one more subset in $\F$, and $\hb(\F)\le
2+\frac{1}{n}\le 2.2$. We will later give examples of families
with Lubell function $>2.2$, so $\F$ cannot satisfy this
condition. Then we may assume by symmetry that $[n]\notin\F$.  We
partition the set $\C$ of full chains into the blocks $\C_{n,i}$,
where the chains $C\in\C_{n,i}$ pass through set
$[n]\setminus\{i\}$. A random full chain in $\C$ is equally likely
to belong to each $\C_{n,i}$, and $\hb(\F)$ is simply the average
over $i$ of the values $\E(|\F\cap C|)$ taken over $C\in\C_{n,i}$,
viewed as the Lubell function for the subsets of
$[n]\setminus\{i\}$. That is, $\hb(\F)$ is the average of $n$
terms, each of which is at most $d_{n-1}$. Hence, $d_n=\hb(\F)\le
d_{n-1}$.


Returning to the calculations, for $d_5$ we note that since
$\hb(\F)$ is a sum of terms, each 1 or $1/5$ or $1/10$, $d_5$ is a
multiple of $1/10$, and hence at most $2.3$ since it is at most
$7/3$.  Then $d_7$, which is similarly a multiple of $1/105$, must
be less than $2.3$, and hence at most $241/105< 2.2953$, and so by
Lemma~\ref{lubbound}, for $n\ge7$, $|\F|< 2.2953\nchn$, which
implies the theorem.  \qed

\medskip
Next is the result for $\D_k$-free families for general $k$.\medskip

\noindent{\bf Proof of Theorem~\ref{dk}}:

Let $n,k\ge2$, and define $m:=\lceil \log_2(k+2)\rceil$.

For the lower bounds, consider $\F=\B(n,m)$. We have $|B-A|\le
m-1$ for any two subsets $A\subseteq B$ in $\F$. There are at most
$2^{m-1}-2$ subsets $S$ satisfying $A\subset S \subset B$. Hence
$\F$ is $\D_k$-free. So $m\le e(\D_k)\le \pi(\D_k)$.

Now we derive the upper bounds. Let $\F$ be a largest $\D_k$-free
family in $\B_n$. We take what we call the {\em min-max
partition\/} of the set $\C$ of full chains in $\B_n$ according to
$\F$:  For subsets $A\subseteq B\subseteq [n]$ with $A,B\in \F$,
the block $\C_{A,B}$ consists of the full chains $C$ such that the
smallest and the largest subsets in $\F\cap C$ are $A$ and $B$,
respectively. We denote by $\C_{\emptyset}$ the block of full
chains that do not meet $\F$ at all.  For $C\in\C_{\emptyset}$, we
have $|\F\cap C|=0$.

We now bound the expected size of $|\F\cap C|$ for a random chain
$C$ in $\C_{A,B}$. If $|B-A|\le m-1$, then this is at most $m$
immediately.

For the remainder, assume $|B-A|\ge m$.  We use the Lubell
function Lemma~\ref{lub} to calculate $\E(|\F\cap C|)$ by adding
the contributions of each subset $S\in \F\cap [A,B]$, which is
$1/\binom{|B-A|}{|S-A|}$. Since $\F$ is $\D_k$-free and
contains both $A$ and $B$, there are at most $k-1$ subsets $S\in
\F\cap [A,B]$ besides $A$ and $B$. Then $\E(|\F\cap C|)$ is
maximized if we take the $k-1$ terms with largest contribution,
{\em i.e.}, with minimum $\binom{|B-A|}{|S-A|}$, which means
the sets $S$ closest to the ends $A$ or $B$, so with $|S-A|$ equal
to 1 or $|B-A|-1$, then 2 or $|B-A|-2$, and so on.  The
contribution from each full level we include is then one.

For the case (1), where $k-1\le 2^{m}-2-{\mchm}$, we see that for
$|B-A|=m$, the $k-1$ terms are at most enough to account for all
subsets $S\in [A,B]$ with $|S-A|$ not equal to $\lfloor
m/2\rfloor$, that is, we get Lubell function at most $m$ (when we
include the terms for $A$ and $B$).  For $|B-A|>m$, since the
levels working up from $A$ or down from $B$ are larger, the $k-1$
terms are no longer sufficient to cover as many full levels, and
the Lubell function is strictly less than $m$.  Since every block
in our partition has expected value at most $m$, we conclude that
$\hb(\F)\le m$. Lemma~\ref{lubbound} gives us $|\F|\le
\Sigma(n,m)$. Furthermore, we also have $e(\D_k)\ge m$. Hence
$m=e(\D_k)=\pi(\D_k)\le \hb(\F)\le m$. By Lemma~\ref{exfam} the
extremal family $\F$ must in fact be $\B(n,m)$.

For the case (2), where $k-1> 2^{m}-2-{\mchm}$, we see that for
$|B-A|=m$, the largest sum of $k-1$ terms leads to Lubell function
at most $m+1-(2^{m}-k-1)/\mchm$.  As in case (1), if $|B-A|>m$,
then since the levels working up from the bottom or down from the
top in $[A,B]$ are larger, the Lubell function is strictly less
than this bound. The bound holds for every block $\C_{A,B}$ of the
min-max partition. Therefore, if $\pi(\D_k)$ exists, $\pi(\D_k)\le
\hb(\F)\le m+1-(2^{m}-k-1)/{\mchm}$. \qed

\medskip
Now we use the min-max partition of the set of full chains to
prove the Harp Theorem.\medskip

\noindent{\bf Proof of Theorem \ref{harp}}:

We argue that $\hb(\F)\le l_1-1$ for any
$\Ha(l_1,\ldots,l_k)$-free $\F$ using induction on $k$. The case
$k=1$ concerns a family $\F$ that contains no chain of height
$l_1$, for which we get immediately that $\hb(\F)=\E(|\F\cap
C|)\le \max_C |\F\cap C|\le l_1-1$ (which implies Erd\H os's
Theorem~\ref{Erd}).

Let $k\ge2$, and assume the bound on $\hb$ for harps with $k-1$
paths.  Let $\F$ be an $\Ha$-free family of subsets of $[n]$,
where $\Ha=\Ha(l_1,\ldots,l_k)$, and consider a block $\C_{A,B}$
in the min-max partition of the set of full chains $\C$ induced by
$\F$. Let $t$ be the largest height of any chain in $\F\cap[A,B]$.
If $t<l_1$ we get that for full chains $C$ in this block,
$\E(|\F\cap C|)\le l_1-1$.

Otherwise, $t\ge l_1$.  Consider a largest chain $Z$ in $\F\cap
(A,B)$, say $S_1\subset \cdots \subset S_{t-2}$, where $A\subset
S_1$ and $S_{t-2}\subset B$. Let $\F'$ be $\F\cap[A,B]$ with the
sets in $Z$ removed.

Then $Z$ and $\F'$ are disjoint and $\E(\F\cap C)$ for random full
chains $C$ in this block is the sum of $\E(Z\cap C)$ and
$\E(\F'\cap C)$.  For the $Z$ term, by Lemma~\ref{lub} we get
$\sum_i 1/\binom{|B-A|}{|S_i-A|}\le (t-2)/|B-A|<1$.  For the
other term, we observe that $\F'$ is $\Ha(l_2,\ldots,l_k)$-free in
the Boolean lattice of subsets of $[A,B]$.  By induction on $k$,
$\hb(\F')\le l_2-1\le l_1-2$.  So $C$ meets $\F'$ on average at
most $l_1-2$ times. Combining terms, we find that $C$ meets $\F$
at most $l_1-1$ times on average for $C$ in this block, and hence
for all random full chains $C\in \C$.  We have that $\hb(\F)\le
l_1-1$.  By Lemma~\ref{lubbound} we get that
$\La(n,\Ha(l_1,\ldots,l_k))\le\Sigma(n,l_1-1)$.

The family $\B(n,l_1-1)$ achieves the upper bound just given,
since it does not contain an $l_1$-chain, and is thus
$\Ha(l_1,\ldots,l_k)$-free.  We see that
\[
e(\Ha(l_1,\ldots,l_k))=\pi(\Ha(l_1,\ldots,l_k))=l_1-1.
\]
Moreover, by Lemma~\ref{exfam}, the only harp-free family of
maximum size is $\B(n,l_1-1)$. \qed


\section{Proof of $\D_2$ Theorem~\ref{d2better}}

We investigate the structure of $\D_2$-free families with maximum
Lubell function, and use this information to improve our earlier
bound.  Before proving Theorem~\ref{d2better}, we continue from
the proof of Proposition~\ref{d2} in the last section, assuming
all notation and facts from that.

We adopt the notation that for any families $\F_1,\ldots,\F_m$ of
sets, $\F_1\vee\cdots\vee\F_m$ denotes the family
$\{F_1\cup\cdots\cup F_m\mid \forall i\, F_i\in \F_i\}$. Given
disjoint sets $S,T$ we define the following three constructions:

\begin{description}
\item{\bf Construction $C_1(S,T)$:} $\F=\{\emptyset\} \cup
\binom{S}{1} \cup \left(\binom{S}{1}\vee\binom{T}{1}\right) \cup \binom{T}{2}$.

\item{\bf Construction $C_2(S,T)$:} $\F=\{\emptyset\} \cup
\binom{S}{2} \cup \binom{T}{2} \cup \left (\binom{S}{2}\vee\binom{T}{1}\right)
\cup \left( \binom{S}{1}\vee\binom{T}{2}\right)$.

\item{\bf Construction $C_3(S,T)$:} $\F=\binom{[n]}{1}\cup
\binom{S}{2} \cup \binom{T}{2} \cup \left( \binom{S}{2}\vee \binom{T}
{1}\right) \cup \left( \binom{S}{1}\vee \binom{T}{2}\right)$.
\end{description}
We will typically partition $[n]$ into subsets $S,T$ in using
these constructions, and we write $C_i(s,n-s)$ for
$C_i([s],[n]\setminus[s])$, for integers $s$, $0<s<n$.  The
families above are $\D_2$-free and each $\hb(C_i(s,n-s))=2+
\frac{s(n-s)}{n(n-1)}$. For $n\ge 2$, the maximum value over $s$
is $2+\frac{\lceil n/2\rceil\lfloor n/2\rfloor}{n(n-1)}> 2.25$,
achieved by $s=\lceil n/2\rceil$ or $\lfloor n/2\rfloor$.

In our approach the key to proving Theorem~\ref{d2better} is to
focus on $\D_2$-free families $\F$ that contain $\emptyset$. Let
$\delta_n$ be the maximum value of $\hb(\F)$ for all such
families.
Definitions give that $2+ \frac{\lceil n/2\rceil\lfloor
n/2\rfloor}{n(n-1)}\le \delta_n\le d_n$. Even though we do not
obtain the values of $d_n$, we can obtain $\delta_n$ for $n$ up to
12.  This technical information (including the extremal families
for $\delta_n$) makes up the following lemma, which is the hard
part in proving the Theorem.

\begin{lemma}\label{l:structure}
The sequence $\{\delta_n\}$ satisfies the following properties.

(1) It is nonincreasing for $n\ge 4$.

(2) For $4\le n\le 12$, if $\F$ contains $\emptyset$ and
$\hb_n(\F)\ge 2\frac{3}{11}$, then up to relabeling elements of
$[n]$, $\F$ is $C_1(s,n-s)$ for $s=\lfloor \frac{n}{2}\rfloor$ or
$\lceil \frac{n}{2}\rceil$, or
$C_2(\lfloor\frac{n}{2}\rfloor,\lceil\frac{n}{2}\rceil)$. Hence,
$\delta_n=2+\frac{\lceil n/2\rceil\lfloor n/2\rfloor}{n(n-1)}$.
\end{lemma}

\noindent{\bf Proof of Lemma \ref{l:structure}:}

To show (1), let $\F$ be a $\D_2$-free family of subsets of $[n]$
such that $\hb_n(\F)=\delta_n$. For $n\ge 5$, if both $\emptyset$
and $[n]$ are in $\F$, then we have only one more subset in $\F$,
and $\hb_n(\F)\le 2+\frac{1}{n}\le 2.2$. Thus $[n]\notin\F$. Then
similar to the proof of Proposition~\ref{d2}, we partition the set
$\C$ of full chains into the blocks $\C_{n,i}$, where the chains
$C\in\C_{n,i}$ pass through set $[n]\setminus\{i\}$.  Again,
$\hb(\F)$ is the average over $i$ of the values $\E(|\F\cap C|)$
taken over $C\in\C_{n,i}$, viewed as the Lubell function for the
subsets of $[n]\setminus\{i\}$. That is, $\hb_n(\F)$ is the
average of $n$ terms, each of which is at most $\delta_{n-1}$. In
other words, for $1\le i\le n$, let $\F_i=\{F\in \F\mid i\not\in
F\}$. Then each $\F_i$ is a $\D_2$-free family in $2^{[n]\setminus
\{i\}}$. We have
\[
\hb_n(\F)=\ave_{1\le i\le n}\left(\ave_{C\in \C_{n,i}}|\F\cap C
|\right)=\frac{1}{n}\sum_{i=1}^n\hb_{n-1}(\F_i).
\]
Hence, $\delta_n=\hb_n(\F)\le \delta_{n-1}$.

We claim the following two facts which are needed in showing (2).

\medskip
\noindent{\bf Claim 1}: The inequality
$2+\frac{\lceil n/2\rceil\lfloor n/2\rfloor}{n(n-1)}-
\frac{1}{\binom{n}{3}}<2\frac{3}{11}$ holds for all $n\le 12$.

This can be verified by a simple computation.

\medskip
\noindent{\bf Claim 2}: For $4\le n\le 12$, suppose $\F\subset C_1(s,n-s)$ with
$\hb_n(\F)\ge2\frac{3}{11}$.  Then $\F=C_1(s,n-s)$ with
$s=\lceil\frac{n}{2}\rceil$ or $\lfloor\frac{n}{2}\rfloor$.
Similarly, if $\F\subset C_2(s,n-s)$ with
$\hb_n(\F)\ge2\frac{3}{11}$, then
$\F=C_2(\lfloor\frac{n}{2}\rfloor,\lceil\frac{n}{2}\rceil)$(The same as
$C_2(\lceil\frac{n}{2}\rceil,\lfloor\frac{n}{2}\rfloor)$ by relabeling the elements).

One can calculate that if $|s-(n-s)|>1$, then
$\hb_n(C_i(s,n-s))<2\frac{3}{11}$. Furthermore, if $\F\subsetneq
\C_1(s,n-s)$ with $s=\lceil\frac{n}{2}\rceil$ or
$\lfloor\frac{n}{2}\rfloor$, then $\hb_n(\F)\le 2+\frac{\lceil
n/2\rceil\lfloor n/2\rfloor}{n(n-1)}-\frac{1}{\binom{n}{2}}
<2\frac{3}{11}$. Similarly, if $\F\subsetneq
C_2(\lfloor\frac{n}{2}\rfloor,\lceil\frac{n}{2}\rceil)$, then
$\hb_n(\F)\le 2+\frac{\lceil n/2\rceil\lfloor n/2\rfloor}{n(n-1)}
-\frac{1}{\binom{n}{3}}<2\frac{3}{11}$. So Claim 2 holds.

\medskip
 We show (2) by induction on $n$.  When $n=4$, it can be
directly verified by enumeration. There are $17$ classes(up
to relabeling of elements of $[n]$) of $\D_2$-free families
containing $\emptyset$. The classes $C_1(2,2)$ and $C_2(2,2)$
satisfy $\hb(\F)=2\frac{1}{3}$ while the rest of them have
$\hb(\F)$ at most $2\frac{1}{4}$ which is less than
$2\frac{3}{11}$.

Assume $n\ge5$ and the statements are true for $n-1$. Now we
consider a $\D_2$-free family $\F\subset 2^{[n]}$ satisfying
$\hb_n(\F)\ge 2\frac{3}{11}$ and $\emptyset \in \F$. Again, the
full set $[n]$ is not in $\F$. Otherwise, $\F$ contains at most
one more subset other than $\emptyset$ and $[n]$, and
$\hb_n(\F)\le 2+\frac{1}{n}<2\frac{3}{11}$. Since $\hb_{n}(\F)\ge
2\frac{3}{11}$, there exists $i$ so that $\hb_{n-1}(\F_i)\ge
2\frac{3}{11}$. We may assume $\hb_{n-1} (\F_n)\ge 2\frac{3}{11}$.
By inductive hypothesis, $\F_n$ is $C_1(\lfloor
\frac{n-1}{2}\rfloor, \lceil \frac{n-1}{2}\rceil)$, $C_1(\lceil
\frac{n-1}{2}\rceil, \lfloor \frac{n-1}{2}\rfloor)$, or
$C_2(\lfloor \frac{n-1}{2}\rfloor, \lceil \frac{n-1}{2}\rceil)$.
We consider two cases.

\medskip
\noindent {\bf Case 1:} $\F_n=C_1(S,T)$ where $|S|=\lfloor
\frac{n-1}{2}\rfloor$ or $\lceil \frac{n-1}{2}\rceil$.

It remains to decide the subsets in $\F\setminus \F_n$. Here are two
subcases depending on whether $\{n\}$ is in $\F$.

{\bf Subcase 1a:} $\{n\}\in \F$.

Since $\F$ is $\D_2$-free, it
contains no subsets of forms $\{s_1,n\}$, $\{s_1,s_2,n,\ldots\}$,
$\{s_1,t_1,n,\ldots\}$, and $\{t_1,t_2,n,\ldots\}$ for $s_i, \in
S$ and $t_i\in T$. Thus, $\F \subset C_1(S\cup \{n\}, T)$. Since
$\hb_n(\F)\ge 2\frac{3}{11}$, we conclude that $|S|+1$ must be
$\lfloor\frac{n}{2}\rfloor$ or $\lceil \frac{n}{2}\rceil$. Thus,
by relabeling elements of $[n]$ we have that $\F=C_1(\lceil
\frac{n}{2}\rceil,\lfloor \frac{n}{2}\rfloor)$.

{\bf Subcase 1b:} $\{n\}\not\in \F$.

Let $S'=\{s\in S\mid
\{s,n\}\in \F\}$ and $T'=\{t'\in T\mid \{t',n\}\in \F\}$. Since
$\F$ is $\D_2$-free, $\F$ cannot have subsets of forms
$\{s_1,s_2,n,\ldots\}$, $\{s_1,t_1,n,\ldots\}$, $\{t_1, t',
n,\ldots\}$, and $\{t_1,t_2, t_3,\ldots\}$ for $s_i\in S$, $t_i\in
T$, and $t'\in T'$. Equivalently,
\begin{center}
$\F\subset \F_n\cup \left( \binom{S'}{1}\vee \{\{n\}\} \right) \cup
\left (\binom{T'}{1}\vee\{\{n\}\} \right) \cup \left(\binom{T\setminus T'}{2}\vee \{\{n\}\}\right)$.
\end{center}
Then
\begin{align*}
  \hb_n(\F)&\le \left(1 + \frac{|S|}{n}+
  \frac{|S||T|+ \binom{|T|}{2}}{\binom{n}{2}}\right)+
  \frac{|S'|+|T'|}{\binom{n}{2}} + \frac{\binom{|T-T'|}{2}}{\binom{n}{3}}\\
&\le  1 + \frac{|S|}{n} +
\frac{|S||T|+ \binom{|T|}{2} + |S'|}{\binom{n}{2}} + f(|T'|).
\end{align*}
Here $f(|T'|)= \frac{\binom{|T-T'|}{2}}{\binom{n}{3}}
+\frac{|T'|}{\binom{n}{2}}$ is a quadratic function of $|T'|$
defined on the integer points of the interval $[0,|T|]$. Its
maximum is reached at one of the two ends, namely $|T'|=0$ or
$|T'|=|T|$.

\noindent{\bf Claim 3}: If $\hb(\F)\ge 2\frac{3}{11}$, then $|T'|=|T|$.

For $n=5$, we have $|S|=|T|=2$. If $|T'|< |T|$, then $f(|T'|)\le
\max\{f(0),f(|T|-1)\}=\frac{1}{10}$. Thus, $\hb(\F)\le 1 +
\frac{|S|}{n} + \frac{|S||T|+ \binom{|T|}{2} + |S'| }{\binom{n}{2}}
 + f(|T'|)<2\frac{3}{11}$, which contradicts our assumption,
and so $|T'|=|T|$.

For $n=6$, we have either $|S|=2$ and $|T|=3$, or else $|S|=3$,
$|T|=2$. If $|T'|<|T|$, then $f(|T'|)\le
\max\{f(0),f(|T|-1)\}=\frac{3}{20}$ for $(|S|,|T|)=(2,3)$, and
$f(|T'|)\le \max\{f(0),f(|T|-1)\}=\frac{1}{15}$ for
$(|S|,|T|)=(3,2)$. By direct computation, both cases give
$\hb(\F)\le 1 + \frac{|S|}{n} + \frac{|S||T|+ \binom{|T|}{2} +
|S'| }{\binom{n}{2}} + f(|T'|)<2\frac{3}{11}$, which is again a
contradiction.

For $7\le n\le 12$, both
$f(0)$ and $f(|T|-1)$ are at most
$\frac{|T|}{\binom{n}{2}}-\frac{1}{\binom{n}{3}}$.
Thus,
\begin{align*}
\hb_n(\F)&\le 1 + \frac{|S|}{n} + \frac{|S||T|+ \binom{|T|}{2} + |S'|
+|T|}{\binom{n}{2}}- \frac{1}{\binom{n}{3}}\\
&\le \hb_n(C_1(|S|,|T|+1))-\frac{1}{\binom{n}{3}}\\
&\le 2+\frac{\lceil n/2\rceil\lfloor n/2\rfloor}{n(n-1)}
-\frac{1}{\binom{n}{3}}< 2\frac{3}{11}.
\end{align*}
This contradiction again proves $|T'|=|T|$, and completes the
proof of Claim 3.

\medskip
Hence, $\binom{T\setminus T'}{2}\vee\{\{n\}\}$ is a null
family. Namely, $\F\subset C_1(S,T\cup\{n\})$. By the condition
$\hb_n(\F)\ge 2\frac{3}{11}$, we conclude by relabeling elements
of $[n]$ that $\F=C_1(\lfloor\frac{n}{2}\rfloor,\lceil
\frac{n}{2}\rceil)$, which is one of the listed possibilities in
(2).

\medskip
\noindent {\bf Case 2:} $\F_n=C_2(S,T)$ where $|S|=\lfloor
\frac{n-1}{2}\rfloor$ and $|T|=\lceil \frac{n-1}{2}\rceil$.

We determine what are the possible subsets in $\F\setminus \F_n$.
Consider the two subcases depending on whether $\{n\}\in\F$.

{\bf Subcase 2a:} $\{n\}\in \F$.

Since $\F$ is $\D_2$-free, $\F$
cannot contain subsets of forms $\{s_1,s_2,n,\ldots\}$,
$\{t_1,t_2,n,\ldots\}$, and $\{u,v,w,n,\ldots\}$ for $s_1, s_2\in
S$, $t_1,t_2\in T$ and $u,v,w \in [n]$. Let $S'=\{s\in S\mid
\{s,n\}\in \F\}$ and $T'=\{t\in T\mid \{t,n\}\in \F\}$. Then
\begin{center}
 $\F \subset \F_n \cup \{\{n\}\} \cup \left(\binom{S'}{1}\vee
  \{\{n\}\}\right) \cup \left(\binom{T'}{1}\vee \{\{n\}\}\right) \cup
  \left(\binom{S\setminus S'}{1}\vee \binom{T \setminus T'}{1}\vee\{\{n\}\}\right)$.
\end{center}
We have
\begin{align*}
\hb_n(\F) &\le \left(1+\frac{\binom{|S|}{2}+ \binom{|T|}{2}}{\binom{n}{2}}
+ \frac{\binom{|T|}{2}|S| + \binom{|S|}{2}|T|}{\binom{n}{3}}\right) + \frac{1}{n}\\
&+ \frac{|S'|+|T'|}{\binom{n}{2}} + \frac{(|S| -|S'|)(|T|-|T'|)}{\binom{n}{3}}\\
&= \hb_n(C_2(|S|+1,|T|))-b,
\end{align*}
 where we see that
 $b:=-\frac{1}{n}+ \frac{|S|-|S'|-|T'|}{\binom{n}{2}}
+\frac{\binom{|T|}{2}+ |S'||T|+|S||T'|-|S'||T'|}{\binom{n}{3}}$ is
a bilinear function of $(|S'|,|T'|)$ defined on $[0,|S|]\times
[0,|T|]$. To find the extremal values of $b$ it suffices to check
the four corner points $(0,0)$, $(|S|,0)$, $(0,|T|)$, $(|S|,
|T|)$. We find the minimum value of $b$, at $(|S|,|T|)$, is
$\frac{|S||T|-\binom{|S|}{2}}{3\binom{n}{3}} \ge
\frac{1}{\binom{n}{3}}$.  Hence,
\[
\hb_n(\F) \le 2+\frac{\lceil n/2\rceil\lfloor n/2\rfloor}{n(n-1)}
- \frac{1}{\binom{n}{3}}<2\frac{3}{11},
\]
which contradicts our assumption, so this subcase
is impossible.

{\bf Subcase 2b:} $\{n\}\not\in \F$.

Similar to subcase 2a, let $S'=\{s\in S\mid \{s,n\}\in \F\}$ and
$T'=\{t\in T\mid \{t,n\}\in \F\}$. The family $\F$ cannot have
subsets of forms $\{s,s',n,\ldots\}$, $\{s',t',n,\ldots\}$, $\{t,
t',n,\ldots\}$, and $\{u,v,w,n,\ldots\}$ for $s\in S$, $t\in T$,
$s'\in S'$, $t'\in T'$ and $u,v,w\in [n]$. Then
\begin{center}
$\F \subset \F_n  \cup \left(\binom{S'}{1}\vee\{\{n\}\}\right) \cup
\left(\binom{T'}{1}\vee\{ \{n\}\}\right) \cup \left(\left(\left(\binom{S}{1}\vee \binom{T}{1}\right)
\setminus \left(\binom{S'}{1}\vee \binom{T'}{1}\right)\right)\vee\{\{n\}\}\right)
\cup\left(\binom{S\setminus S'}{2}\vee \{\{n\}\}\right) \cup
\left(\binom{T\setminus T'}{2}\vee \{\{n\}\}\right)$.
\end{center}
We have
\begin{align*}
\hb_n(\F)&\le \left(1 + \frac{\binom{|S|}{2}+ \binom{|T|}{2}}{\binom{n}{2}}
+\frac{\binom{|T|}{2}|S| + \binom{|S|}{2}|T|}{\binom{n}{3}}\right)
+\frac{|S'| +|T'|}{\binom{n}{2}}\\
& +\frac{|S||T| -|S'||T'| + \binom{|S\setminus S'|}{2} +
\binom{|T\setminus T'|}{2}}{\binom{n}{3}}\\
&= \hb_n(C_2(|S|+1,|T|))- \frac{g}{\binom{n}{3}},\mbox{ where}\\
g&= \binom{n}{3} \left(\frac{|S|-|S'|-|T'|}{\binom{n}{2}}\right)+\binom{|T|}{2}
+ |S'||T'| -\binom{|S\setminus S'|}{2} - \binom{|T\setminus T'|}{2}\\
&= \varepsilon_1 |T'|+ \varepsilon_2|S\setminus S'|,\mbox{ with}
\end{align*}
\[
\varepsilon_1=\frac{2(n-2)}{3}-\frac{|T'|-1}{2}
-\frac{5(|S\setminus S'|)}{6}\mbox{ and }
\varepsilon_2=\frac{n-2}{3}-\frac{|S|-|S'|-1}{2}-\frac{|T'|}{6}.
\]

If $(|S'|,|T'|)=(|S|,0)$ or $(0,|T|)$, then $\F \subset
C_2(S\cup\{n\},T)$ or $\F\subset C_2(S,T\cup\{n\})$. Since
$\hb_n(\F) \ge 2\frac{3}{11}$, we have, by relabeling elements of
$[n]$, $\F=C_2(\lfloor \frac{n}{2}\rfloor,\lceil
\frac{n}{2}\rceil)$, which is another of the alternatives listed
in (2).

\medskip
\noindent{\bf Claim 4}: When $(|S'|,|T'|)\neq(|S|,0)$ or $(0,|T|)$,
then $g\ge 1$.

Recall that $0\le |S'|\le |S|=\lfloor\frac{n-1}{2}\rfloor$ and
$0\le |T'|\le |T|=\lceil\frac{n-1}{2}\rceil$.

Suppose $n=2k$, so we have $|S|=k-1$ and $|T|=k$, $k\ge 2$.
Rewrite

\[
\varepsilon_1=(3|T\setminus T'| + 5|S'|)/6\ \mbox{ and }
\varepsilon_2 =(|T\setminus T'| + 3|S'| +2)/6.
\]
Note that $|S'|$ and $|T\setminus T'|$ are not both zero, nor are
$|T'|$ and $|S\setminus S'|$ both zero.  One can see that
$\varepsilon_1$ and $\varepsilon_2$ are each at least $1/2$, and
so $g\ge1$ unless $|T'|=1$ and $|S\setminus S'|=0$,  or $|T'|=0$
and  $|S\setminus S'|=1$. But either pair of conditions increases
the $\varepsilon$'s and still leads to $g\ge 1$.

Else, suppose $n=2k+1$, $|S|=|T|=k$ and $k\ge 2$. This time
rewrite
\[
\varepsilon_1 =(3|T\setminus T'| +5|S'|-1)/6\ \mbox{ and }
\varepsilon_2 =(|T\setminus T'|+3|S'|+1)/6.
\]
Again it is simple to check that $g\ge1$, and Claim 4 holds.

From Claim 4, we have if $(|S'|,|T'|)\neq(|S|,0)$ or $(0,|T|)$,
then once again we get the contradiction
\[
\hb_n(\F)\le 2+\frac{\lceil n/2\rceil\lfloor n/2\rfloor}{n(n-1)}
-\frac{1}{\binom{n}{3}}< 2\frac{3}{11}.
\]
This completes the Case 2 and the proof of the Lemma. \qed

Now we are ready to prove our improved bound.

\medskip

\noindent{\bf Proof of Theorem~\ref{d2better}}:
Let $\F$ be a
$\D_2$-free family of subsets of $[n]$. Partition $\F$ into $\F_k$
and $\F'$ such that $\F_k$ contains subsets of sizes in $[k,n-k]$
where $k=\frac{n}{2}- 2\sqrt{n\ln n}$, and $\F'=\F\setminus \F_k$.
We know that $|\F'|/\nchn<\frac{2}{n^2}$ for large $n$
(see~\cite{GriLu}, Lemma 1).

Now concentrate on the family of sets near the middle, $\F_k$.  We
take what we call the {\em min partition\/} of the set $\C$ of
full chains in $\B_n$:  Let $\C_{\emptyset}$ be the block
containing the full chains that do not meet $\F_k$ at all.  For
each subset $A\in \F_k$, let $\C_A$ be the block containing all
full chains $C$ having $A$ as the minimal element in $\F_k\cap C$.
We see that the average number of times a chain in $\C_A$ meets
$\F_k$ is obtained by considering only the subsets in $\F_k$ that
contain $A$ and viewing them (after removing $A$ from each) as a
diamond-free family of subsets of $[n]\setminus A$ containing $\emptyset$.
We deduce that $\ave_{C\in\C_A}(|\F_k\cap C|)\le\delta_{n-|A|}$.
For $n$ large enough, we have $n-|A|>\frac{n}{2}- 2\sqrt{n\ln
n}\ge 12$. Since the $\delta_n$ are nonincreasing, we have for
large $n$ that $\hb(\F_k)\le \delta_{12}=2\frac{3}{11}$. It
follows that for sufficiently large $n$, all $\D_2$-free families
$\F$ in $\B_n$ satisfy
\[
\frac{|\F|}{\nchn}=\frac{|\F_k|+|\F'|}{\nchn}\le 2\frac{3}{11}+
\frac{2}{n^2}.
\]
Consequently, if it exists, the limit
$\pi(\D_2)\le 2\frac{3}{11}$.\qed

\section{Further Research}

Beyond diamonds $\D_k$, we continue to investigate why the limit
$\pi(P)$ exists for general posets $P$. The methods introduced in
this paper have proven to be useful for determining $\pi(P)$ for
several other small posets $P$, which we are collecting
separately~\cite{GriLi,GriLi2,WTL}. One example is the subposet
$\J$ of $\D_2$ consisting of four elements $A,B,C,D$ with $B<A$ and
$B<C<D$. Forbidding $\J$ is more restrictive than forbidding
$\D_2$.  We show that for $n\ge1$, $\La(n,\J)=\Sigma(n,2)$, and
hence $\pi(\J)=2$.  All known values of $\pi(P)$ satisfy
Conjecture~\ref{conj}.

In order to resolve the asymptotics for diamond-free posets, and
show that $\pi(\D_2)=2$ as expected, it is not enough to work with
the Lubell function due to families such as those in the
constructions $C_i$ described above. These examples show that the
terms in the sequence $\delta_n$, which was shown to be
nonincreasing for $n\ge2$, are at least $2.25$ for all $n\ge2$. We
suspect that the limit $\lim_{n\rightarrow\infty} d_n$, which is
known to exist, is $2.25$.  This would follow from the conjecture
below. Here, the {\em conjugate} of $\F$ is the family
$\overline{\F}=\{\bar F\mid F\in \F\}$ where $\bar F=[n]\setminus
F$ is the complement of $F$.

\begin{conjecture}
  For every $n\ge 4$, the value $\hb(\F)$ of any $\D_2$-free
  family $\F\subset 2^{[n]}$ satisfies $\hb_n(\F)\le2+
  \frac{\lfloor\frac{n^2}{4}\rfloor}{n(n-1)}$ and equality holds if
  and only if, up to relabelling elements of $[n]$, $\F$ or
  $\overline{\F}$ is $C_i(S,T)$ ($i=1,2,3$) with $||S|-|T||\le 1$.
\end{conjecture}

Then how might we reduce our upper bound on $\pi(\D_2)$ (if it
exists) to below $2.25$? The examples above are nowhere near as
large as $2\nchn$, yet have very small sets that make large
contributions to the Lubell function.  To build large diamond-free
families, we should restrict our attention to families with no
small nor large sets, say $\F\subseteq \B(n,k)$ with $k=n-f(n)<n$.
If we could show that $|\F|$ is  at most $(2+o_n(1))\nchn$, for
suitable $f(n)$, we would have $\pi(\D_2)=2$ as we expect, since
most of the $2^n$ subsets are concentrated near the middle rank.

Another indication of the challenge facing us for $\D_2$ is that
we have constructed three $\D_2$-free families for $n=6$ of size
$36$, which is one more than $\Sigma(6,2)$.  This is in contrast
to the values $k$ for which Theorem~\ref{dk} determines
$\La(n,\D_k)$ completely, and its value is exactly $\Sigma(n,m)$.
Thus, the solution for $\D_2$, and probably also $\D_k$ for the
unsettled values $k$, is likely going to be more complicated.


\comments{
\section{Appendix}
\begin{lemma}
Suppose $X$ is a random variable which takes on nonnegative
integer values. Let $f(x)$ and $g(x)$ be two increasing functions.
Then}
\end{document}